\documentclass{article}
\usepackage{geometry}
\usepackage{amsmath,float,amssymb,amsfonts}
\usepackage{url,cite}
\usepackage{multirow}
\usepackage{longtable}
\usepackage{color}
\usepackage{verbatim}
\usepackage[linesnumbered,ruled,lined]{algorithm2e}
\usepackage{amsthm}
\usepackage{mathtools}

\usepackage{pgfplots}
\pgfplotsset{%
,compat=1.12
,every axis x label/.style={at={(current axis.right of origin)},anchor=north west}
,every axis y label/.style={at={(current axis.above origin)},anchor=north east}
}

\tikzset{mynode/.style={inner sep=2pt,fill,outer sep=0,circle}}
\usepackage{enumerate}

\newtheorem{theorem}{Theorem}

\newtheorem{definition}{Definition}

\usepackage{graphicx,tikz,pgf}
\usepackage{enumerate}
\usepackage{authblk}
\usepackage{tikz}
\usepackage{pgfplots}
\usetikzlibrary{positioning}

\date{}

\begin{document}

\title{On the Rudnick and Sarnak's Zeros of principal L-functions and Random Matrix Theory}

\author{Madhuparna Das}
\maketitle

\begin{abstract}
In this srticle we have surveyed the result of Ze\'ev Rudnick and Peter Sarnak on the Zeros of principal L-function and Random Matrix Theory.
\end{abstract}

\section{Introduction}

In this article, we will survey the results of Rudnick and Sarnak~\cite{bib1} on the ``Zeros of Principle L-function and Random Matrix Theory". The main concern of this paper is to study about the zeros of Riemann zeta function as well as more general L-functions\footnote{The Riemann zeta function $\zeta(s)$ is an especial case of L-functions.}. First we see the motivation of studying zeta function as well as L-functions. L-functions are certain analytic functions in number theory which helps with the arithmetic information to study its analytic behaviour.
Let us look at the Riemann zeta function first then we can extend the idea for more general L-functions. Riemann zeta function $\zeta(s)$ is defined by

\begin{align*}
\zeta(s) = \sum_{n=1}^{\infty} \frac{1}{n^s} = \prod_{p} \frac{1}{1-p^{-s}},
\end{align*}
where the product runs over the primes, known as Euler product. The above series converges for $Re(s)>1$. Analytically continued to a meromorphic function on the complex plane $\mathbb{C}$ with a simple pole at $s=1$ and satisfies the functional equation relating $\zeta(s)$ and $\zeta(1-s)$. Now we can move to the L-function; before going into more general L-functions we see the simplest one which is called the Dirichlet L-function and defined by,

\begin{align*}
L(s,\chi) = \sum_{n=1}^{\infty} \frac{\chi(n)}{n^s} = \prod_{p} \frac{1}{1-\chi(p)p^{-s}},
\end{align*}
with $Re(s)>1$ and $\chi$ is being a Dirichlet character modulo $N$. Observing the above equation, we can see that if we take $\chi=1$ then the function $L(s,\chi)$ leads to the $\zeta$ function and we already know its analytical behaviour. For $\chi\ne1$, $L(s,\chi)$ is an entire function of $\mathbb{C}$ with the functional equation which relates $L(s,\chi)$ and $L(1-s, \bar{\chi})$.

The motivation behind the study of zeta function as well as L-function because they provide arithmetic information and due to their analytical behaviour. For example, zeta has a simple pole at $s=1$ which implies that there exist infinitely many primes and the zeros of zeta function indicate the distribution of prime numbers. Similarly, for $\chi\ne1$, $L(s,\chi)$ has no pole and it is nonzero at $s=1$, which indicates that there are infinitely many primes in an arithmetic progression with gcd 1.

\medskip

Let us now look at the algebraic behaviour of Dirichlet L-function, which provides us a new way to look at the L-functions and also, provides information about the Langlands program. As we know the Dirichlet character $\chi$ is an algebraic object. Class field theory asserts that there is a correspondence between primitive Dirichlet characters and Hecke Characters of finite order. The Hecke character of finite order is defined by

\begin{definition}{\bf Hecke Character:}
A Hecke charatcer of finite order ${\chi}_1= \otimes{{\chi}_1}_v$ is a character ${\chi}_1: \mathbb{A}^{\times}/\mathbb{Q}^{\times}\rightarrow\mathbb{C}^{\times}$ whose kernel has finite index in the idele class group $\mathbb{A}^{\times}/\mathbb{Q}^{\times}$.
\end{definition}

Then one can define the Hecke L-function by

\begin{align*}
L(s,{\chi}_1) = \prod_{p} L(s,{{\chi}_1}_p)
\end{align*}
with $Re(s)>1$ and $L(s,{{\chi}_1}_p)$ is Dirichlet L-function if ${{\chi}_1}_p$ is unramified and 1, if ${{\chi}_1}_p$ is ramified (means $p|N$). So the correspondence between Dirichlet and Hecke L-function is given by
\begin{align*}
L(s,\chi) = L(s, {{\chi}_1}).
\end{align*}

There are other L-functions except for Hecke and Dirichlet L-functions. If we look at the Dedekind zeta function it leads us to the Artin L-functions which generalize the Dirichlet L-function. Now the question arieses that is it possible to find a correspondence between generalized Hecke characters and Artin L-function. Well, Laglads conjecture assures that there is a correspondence for which we should look at the automorphic representation of $\pi$ of $GL_n$\footnote{where $n$ is the dimensions of Artin representation $\rho$} over a number field. The automorphic representations are analytic objects so relatively it's easy to study the properties of their L-functions. Indeed, the analogue of Artin’s conjecture is known for automorphic representations, so if one could show Artin representations correspond to automorphic representations (in the sense that their L-functions agree), one could deduce Artin's conjecture. The theory of principle L-functions was developed by Godment and Jacquet~\cite{bib2} for $n\geq3$ (from the case $\mathbb{A}^{\times}$, where $\mathbb{A}^{\times}$ is the adeles of $\mathbb{Q}$)and by Hecke and Mass for $n=2$.

The main interest of the result given by Rudnick and Sarnak was to study the fine structure of the distribution of the non-trivial zeros of such primitive $L(s,\pi)$. If $\rho^{\pi}=\frac{1}{2}+i\gamma^{(\pi)}$ denotes the zeros of such non trivial primitive L-functions $L(s,\pi)$,then by assuming Riemann hypothesis for $L(s,\pi)$ i.e., $\gamma\in{\bf R}$, they have ordered ${\gamma}^{(\pi)}$'s (with multiplicities)

\begin{align*}
\cdots \leq {\gamma^{(\pi)}}_{-2} \leq{{\gamma}^{(\pi)}}_{-1}<0\leq {{\gamma}^{(\pi)}}_{1}\leq{{\gamma}^{(\pi)}}_{2}\cdots
\end{align*}

The main problem of Rudnick and Sarnak's paper is to understand the asymptotic behaviour of $\gamma$ in an interval and their statistical distribution. In the case of Riemann zeta function the calculation given by Montgomery~\cite{bib3} and Odlyzko~\cite{bib4}\cite{bib5} says that the consecutive spacings follow the Gaussian Unitary Ensemble (GUE) distribution from the Random Matrix Theory. In the later sections, we will discuss these statistical distributions and their relation with Random Matrix Theory.

\section{The Correlation Sum $R_n(B_N,f)$}

The Gaussian Unitary Enssemble (GUE) is if ${\delta}_n={\gamma}_{n+1}-{\gamma}_{n}$ are normalized spacings, then for any nice function on $(0,\infty)$ we can expect
\begin{align*}
\frac{1}{N}\sum_{n\leq N} f({\delta}_n) \to \int_{0}^{\infty} f(s)P(s) ds,
\end{align*}
where $P(s)$ is the distribution of consecutive spacing of the eigenvalues of a large random Harmition matrix.

The distribution says,
\begin{align*}
P(s)= \frac{d^2E}{ds^2}(s) \text{where $E(s)$ = det$(I-Q_s)$},
\end{align*}
where $Q_s$ is the trace class operator on $L^2(-1,1)$ with kernel $Q_s(\xi,\eta)=\frac{sin\pi s(\xi-\eta)}{\pi s(\xi-\eta)}$.

The main result given by Rudnick and Sarnak is the computation of the general $n$-level correlation function for the zeros of a primitive principal L-function. They have shown that the answer is universal and is precisely the one predicted by Dyson's computations for the GUE model~\cite{bib6}.

\medskip

The deifiniton of $n$-level correlations says, suppose we have a set $B_N$ of $N$ numbers,
\begin{equation*}
B_N=\{{\gamma}_i: i\in\{1,\ldots,N\}\}
\end{equation*}
holds the ineqaulity,
\begin{align*}
\tilde{{\gamma}}_1\leq \ldots \leq \tilde{{\gamma}}_N.
\end{align*}
The $n$-level correlation function measures the correlation between differences $n$ elements of $B_N$.

\medskip

So this number is
\[
N(a,b)=\#\{k:\tilde{\gamma}_{k+1}-\tilde{\gamma}_k\in[a,b]\}
\]
for any interval $[a,b]$.\\
However it is hard to calculate $N(a,b)$ directly, because it is not easy to tell if two elements are consecutive or not unless we know all the numbers in the sequence. We can put
\[
N'=N'(a,b)=\#\{(\tilde{\gamma},\tilde{\gamma}')\in B_N^2|\tilde{\gamma}<\tilde{\gamma}',\tilde{\gamma}'-\tilde{\gamma}\in [a,b]\}
\]
\[
N''=N''(a,b)=\#\{(\tilde{\gamma},\tilde{\gamma}',\tilde{\gamma}'')\in B_N^3|\tilde{\gamma}<\tilde{\gamma}'<\tilde{\gamma}'',\tilde{\gamma}''-\tilde{\gamma}\in [a,b]\}
\]
and so on, where $N^{(n)}$ is the number of $n$ tuples whose differences between the biggest and the smallest elements belong to $[a,b]$. Using induction, we can write
\[
N(a,b)=N'-N''+....
\]
Since the interval is $[a,b]$, the above sum is finite. It is sufficient to obtain the numbers $N^{(n)}$ for $n=2,3,4,...$.
\begin{definition}{\bf $n$ level correlation}
Let $f$ be a symmetric $n$ variable function given by $f(S)=f(a_1,\ldots,a_n)$ if $S=\{a_1,\ldots,a_n\}$. For the box $Q\in {{\bf R}^{n-1}}$,

\begin{align*}
R_n(B_N,Q) =\frac{1}{N} \#\{j_1,\cdots j_n\leq N distinct : (\tilde{\gamma}_{j_1}-\tilde{\gamma}_{j_2},\ldots,\tilde{\gamma}_{j_{n-1}}-\tilde{\gamma}_{j_n})\in Q\}
\end{align*}

gives,
\begin{align*}
R_n(B_N,f)= \frac{n!}{N}\sum_{\substack{S\subset B_N \\ |S|=n}} f(s)
\end{align*}
\end{definition}

So the function $f$ satisfies the following three conditions
\begin{enumerate}
\item $f(x_1,...,x_n)$ is symmetric;
\item $f(x+t(1,...,1))=f(x)$ for $t\in \mathbf{R}$;
\item $f(x)\to 0$ rapidly as $|x|\to \infty$ in the hyperplane $\sum_j x_j=0$
\end{enumerate}
Knowing the asymptotic behaviour of $R_n(B_N,f)$ as $N\to\infty$ is equivalent to knowing that of the smoothed correlations
\[
R_n(T,f,h)=\sum_{j_1,...,j_n}h\Big(\frac{\gamma_{j_1}}{T}\Big)...h\Big(\frac{\gamma_{j_n}}{T}\Big)f\Big(\frac{L}{2\pi}\gamma_{j_1},...,\frac{L}{2\pi}\gamma_{j_n}\Big)
\]
for a sufficiently rich family of localized cut off functions $h$. Here $L=m\log T$ and $\tilde{\gamma}_j$ is normalized by $\tilde{\gamma}_j=\frac{\gamma_j\log \gamma_j}{2\pi}$ given by Riemann.

\medskip

Dyson determined the density of the limiting $n$-correlation sum $W_n(x_1,...,x_n)$ for GUE model.\\
$W_n(x)$ is a density satisfying $0\leq W_n(x)\leq 1$
\begin{equation*}
W_n(x) = \begin{cases}
0, & \text{if } x_i=x_j \text{ for } i\neq j\\
1, & \text{iff $x_i-x_j\in \mathbb{Z}$ and $x_i\neq x_j$ for all $i\neq j$}
\end{cases}
\end{equation*}
Then $W_n(x_1,...,x_n)=det(k(x_i-x_j))$, where $k(x)=\frac{\sin\pi x}{\pi x}$

\begin{theorem}[Rudnick-Sarnak]\label{thm1}
Let $\pi$ be a cuspidal automorphic representation of $GL_n/\mathbf{Q}$. Assume $m\leq 3$ or the hypothesis\footnote{The hypothesis asserts that for any $k\geq2$, $\sum_{p}\frac{|a_{\pi}(p^k)\log p|^2}{p^k}<\infty$}. Let $f$ satisfy the all three conditions and in addition assume that $\hat{f}(\xi)$ is supported in $\sum_j|\xi_j|<2/m$. Let $g\in C_c^\infty(\mathbf{R})$ and $h(r)=\int_{-\infty}^\infty g(u)e^{iru}du$ (so that $h$ and $f$ are entire). Then as $T\to\infty$
\[
R_n(T,f,h)\sim \frac{m}{2\pi}T \log T\int_{-\infty}^\infty h(r)^n dr\int_{\mathbb{R}^n}f(x)W_n(x)\delta\Big(\frac{x_1+...+x_n}{n}\Big)dx_1...dx_n
\]
Where $\delta(x)$ is the Dirac mass at $0$.
\end{theorem}

Before going into the proof of Theorem~\ref{thm1} let us discuss some more about the principal L-functions and Rankin-Selberg Convolution.

\subsection{Principal L-function}

Rudnick and Sarnak focused especially on the automorphic L-functions on $GL_m$ for $m\geq3$. For the lower ranks, L-functions (e.g. Dirichlet L-function) are classical. It can be shown that Dirichlet L-functions satisfy the anlytic properties using integeral representation which plays an important role in the theory of L-functions. If $\chi$ is a Dirichlet character we can write
\[
L(s,\chi)=\int_{0}^{\infty} {\Phi}_{\alpha} t^{s/2} \frac{dt}{t}
\]
for a function ${\Phi}_{\alpha}$.

One can break the integral and by using the change of variables one can obtain the functional equation for $L(s,\chi)$ relating with $L(s,\chi-1)$.

\medskip

Let $\pi$ be a cuspidal automorphic representation of $GL_m/{\bf Q}$. Consider the factor $\pi={\otimes}_p\pi_p$ where $\pi$ is a smooth irreducible (infinite-dimensional) representation of $GL_m/{\bf Q}$ for each $p$. Then we can define $L(s,\pi)$ associated with Euler product,
\[
L(s,\pi) = \prod_{p} L(s,{\pi}_p)
\]
which is the product of local factors and for almost all $p$, $\pi$ is an unramified principal series. That means there is such a conjugacy class which is parametrized by its eigenvalues ${\alpha}_{\pi}(j,p)$ (for $j=1,\ldots,m$).
The local factors $L(s,{\pi}_p)$ for the unramified primes are given by,
\[
L(s,{\pi}_p) = det(I-p^{-s}A_{\pi}(p))^{-1}=\prod_{j=1}^m (1-{\alpha}_{\pi}(j,p)p^{-s})^{-1}
\]

So, $L(s,\pi_p)$ should be the reciprocal of a polynomial of degree $m$ in $p^{-s}$ when $\pi_p$ is unramified. In general, at any prime $p$, $L(s,\pi_p)$ should be the reciprocal of a polynomial of degree less than equal to $m$ in $p^{-s}$, so we will say $L(s,\pi)$ is an L-function of degree $m$. But the case of $L(s,\pi_{\infty})$, it is a product of Gamma functions (or local archimedean factor) given by
\[
L(s,\pi_{\infty}) = \prod_{j=1}^m \Gamma_{\bf R}(s+\mu_{\pi}(j))
\]

where $\Gamma_{\bf R}(s)={\pi}^{-s/2}\Gamma(\frac{s}{2})$ and $\{\mu_{\pi}(j)\}$ is a set of $m$ numbers associated to $\pi_{\infty}$.

With all the local factors defined one can turn to the functional equation given by Godement and Jacquet which can be deduced from the Langlands–Shahidi method.

\medskip

Now we can write that there is pricipal L-function $L(s,\pi)$ which is associated with $\pi$ and entire, also satisfies the functional equation
\[
\Phi(s,\pi)=\epsilon(s,\pi)\Phi(1-s,\tilde{\pi}) \\
\epsilon(s,\pi)=\tau(\pi)Q_{\pi}^{-s}
\]
where $Q_{\pi}>0$ is the conductor of $\pi$ and $\Phi(s,\pi)$ satisfies
\[
\Phi(s,\pi)=L(s,\pi_{\infty})L(s,\pi)
\]

Godement and Jacquet had not computed the local factors of $L(s,\pi_p)$ for all the cases when $\pi_p$ is ramified. In particular, the non-archimedean case which is not really needed for this article. Note that the non-trivisl zeros of $L(s,\tilde{\pi})$ are realted to those of $L(s,\pi)$ via $s\to1-s$. According to the Riemann Hypothesis, $Re(\rho_{\pi})=1/2$, then we have the counting function
\[
N_{\pi}(T)=\# \{\rho_{\pi}:|Im\rho_{\pi}|<T\}
\]
is asymptotic to $\frac{m}{\pi}T\log T$.

We give a brief overview of the computation. Let $N'_{\pi}(T)$ be the number of zeros of $\rho_{\pi}$, then we can write
\[
N_{\pi}(T) = N'_{\pi}(T)+N'_{\tilde{\pi}}(T)+ \text{Error}
\]

We can assume that $\Lambda_{\pi}(T)$ does not vanish on $Im(s)=T$ and we know that $-T<|Im\rho_{\pi}|<T$. From this we can get a logarithmic derivative of $\frac{\Lambda'}{\Lambda}(T)$. Now integrating over the rectangle and computing using contour integral, we can obtain the terms involving gamma factors. From the asymptotic behaviour of gamma (one can deduce from the Sterling formula), we can obtain that $N_{\pi}(T)\sim\frac{m}{\pi}T\log T$. The term $m/\pi$ is coming due to the degree of the L-function.

\subsection{Rankin-Selberg Convolution}

%We are going to give a brief overview of the Rankin-Selberg convolution.

%\section{Rankin-Selberg Convolution}

In this section, we are going to give a brief overview of two famous mathematician's work named by Rankin-Selberg Convolution by Rankin and Selberg. We start from a very basic thing from Calculus, a contribution of Euler, \\

{\bf Euler:} Euler Integral:
\begin{align*}
&\int_{0}^{1} x^{\alpha-1}(1-x)^{\beta-1}dx \\
&= \frac{\Gamma(\alpha)\Gamma(\beta)}{\Gamma(\alpha+\beta)}.
\end{align*}
\medskip
where $\Gamma(\alpha)=\int_{0}^{\infty}e^{-x}x^{\alpha}\frac{dx}{x}$.

\medskip

Now we move our focus on the Selberg integral, in which he found a much deeper result from Euler and it states that

\begin{align*}
& \int_{0}^{1} \ldots \int_{0}^{1} (x_1x_2\ldots x_n)^{\alpha-1} \left((1-x_1)(1-x_2)\ldots(1-x_n)^{\beta-1}\right)|\Delta(x)|^{2\gamma} dx_1\ldots dx_n \\ \\
& =\prod_{j=1}^{0} \frac{\Gamma(1+\gamma+j\gamma)\Gamma(\alpha+j\gamma)\Gamma(\beta+j\gamma)}{\Gamma(1+\gamma)\Gamma(\alpha+\beta+(n+j-1)\gamma)}
\end{align*}
\\
where the term $\Delta(x)=\prod_{1<j\leq i\leq n} (x_j-x_i)$.

{\bf Proof Idea:} For an integer $\gamma\geq0$ we can expand the multivariable integral $\left(\Lambda(x)\right)^{2\gamma}$ and then apply Euler's integral and some critical argument on symmetry. Then collect all $\gamma$ using functions theory which proves the above equation.

It implies Dyson-Mehta Conjecture which states that

\begin{equation*}
\begin{split}
\int_{-\infty}^{\infty}\ldots \int_{-\infty}^{\infty} e^{-\sum_{j=1}^{\infty}|\Delta{(x)}|^{2\gamma}}dx_1\ldots dx_n \\
=(2\pi)^{\alpha/2} \prod_{j=1}^{n} \frac{\Gamma(1+j\gamma)}{\Gamma(1+\gamma)}
\end{split}
\end{equation*}

Later in 1982, Macdonald generalized the Dyson-Mehta conjecture by replacing $\Delta(x)$ with $P(x)$, where $P(x)$ is the product of the distances to the hyperplane of a finite representations group $\mathbb{R}^n$.

\medskip

Gauss sums involve gamma function and Jacobi sum involves $\beta$ function. In 1990 Anderson worked on the Selberg's sums involving Selberg integral.

\medskip

Well Rankin-Selberg L-function was followed by Ramanujan and Hecke

\begin{align*}
&\Delta(q)=q\prod_{n=1}^{\infty}(1-q)^{24} \\
&= q-24q+252q^2-1472q^4\ldots \\
&= \sum_{m=1}^{\infty} \tau(m)q^m
\end{align*}
where $\tau(m)$ is the Ramanujan's Tau function. Ramanujan Tau function follows certain properties: 
\begin{itemize}
\item $\tau(mn)=\tau(m)\tau(n)$ if $gcd(m,n)=1$.
\item $|\tau(p)|=2p^{s/2}$.
\end{itemize}
The above properties of Tau function is also known as ``Ramanujan's Conjecture".

In the case of classical Modular forms, if $q=e^{2\pi iz}$ and $Im(z)>0$, then

\begin{equation*}
\Delta\left(\frac{az+b}{cz+d}\right) = (cz+d)^{12}\Delta(z)
\end{equation*}

for the matix
\[
\begin{bmatrix}
a & b \\
c & d
\end{bmatrix}
\]
$\Delta$ is a modular form of weight 12 for $SL_2(\mathbb{Z})$.

The contribution of Hecke says that

\begin{equation*}
L(s,\tau)=\sum_{n=1}^{\infty} \frac{\tau(n)}{n^{1/2}}n^{-s}\\
=\prod_{p} \left(1-\frac{\tau(p)}{p^{s/2}}p^{-s}+p^{-2s}\right)^{-1}
\end{equation*}
then $L(s,\tau)$ satisfies properties similar to Ramanujan's tau function and its analytic to continuation functional equation. \\

{\bf Rankin-Selberg:}
\begin{align*}
L(s, \tau\times\tau) = \sum_{n=1}^{\infty} \frac{\tau(n)^2}{n^{11}}n^{-1}
\end{align*}
has similar properties. \\

{\bf Idea of Proof:} Let us look at the integral:

\begin{align*}
I(s)= \int_{SL_2(x)^{z-1}} |\Delta(z)|^{2}y^2 E(z,s) \frac{dxdy}{y^2}
\end{align*}

where $E(z,s)=\sum_{(m,n)\ne (0,0)} \frac{y^s}{(mz+n)^{2s}}$.
If we unfold the Rankin-Selberg integral
\[
I(s) = *\times L(s,\tau\times\tau)
\]
where $*$ is some known factors. We understand the integral using Eisenstien Series.

The properties we can conclude are as followes:
\begin{itemize}
\item The trivial upper bound of $\tau(p)=2p^{11/2+\epsilon}$ (an approximation of Ramanujan's tau function). 
\item Langland's put forth a general family of L-functions in association with $\tau$. Hecke and Rankin Selberg are the first two in this family. The conjectured properties of these would imply Ramanujan's Conjecture.
\item The Rankin-Selberg method has been developed more generally in the context of automorphic forms for the group of $n\times n$ matrices. $GL_2$ is the classical case and $GL_n$ is a context test in the theory Jean-Perre-Serre.
\end{itemize}
Now we move our focus to the main result of Rudnick-Sarnak.

\section{Asymptotic behaviour of the correlation sum $R_n(f,T)$}

In this section, we wish to study the asymptotic behaviour of $n$ correlation function

\begin{align}\label{eq1}
R_n(f,T)=\frac{1}{N} \sum_{i_1,\ldots, i_n\leq N}^{*} f(\tilde{\gamma_{i_1}},\ldots \tilde{\gamma_{i_1}})
\end{align}
where $N=N(T)$ and $\sum^{*}$ means sums over distinct indices $i_j$.

It's hard to know about the asymptotic behaviour of $R_n(f,T)$ directly. So, instead of a direct calculation, Rudnick and Sarnak broke the sum into some steps which make the calculation easy. As a first step, let us look at the sum
\begin{align*}
C_n(f,T)=\sum_{i_1,\ldots i_n\leq N} f(\tilde{\gamma}_{i_1},\ldots, \tilde{\gamma}_{i_1})
\end{align*}

Observe that the sum is no longer distinct ordered zeros as in the definition of the correlation function. Rudnick-Sarnak has computed the sum $R_n(f, T)$ by using combinatorial sieve. It's now clear that first, we want to calculate the asymptotic behaviour of the $C_n(f, T)$ sum and later the $R_n(f, T)$ sum. Adding some more steps will make the computation a bit easy. To get the sum $C_n(f,T)$ consider the smooth sum

\begin{align*}
C_n(f,h,T) = \sum_{\gamma_1,\ldots,\gamma_n} h_1\left(\frac{\gamma_1}{T}\right) \ldots h_1\left(\frac{\gamma_n}{T}\right) f\left(\frac{L}{2\pi}\gamma_1,\ldots,\frac{L}{2\pi}\gamma_n\right)
\end{align*}
Set $L=m\log T$ and $h_j(r)$ is a smooth cutoff. Set
\begin{align}\label{eq2}
h_j(r)=\int_{-\infty}^{\infty} g_j(u)e^{iru}du
\end{align}
with $g_j\in C_c^{\infty}({\bf R})$. Now we state and prove the theorem to deduce the asymptotic behaviour of the sum $C_n(f,h,T)$.

\begin{theorem}\label{thm2}
Let $\Phi\in C^1$ be supported in $\sum_{j=1}^{n}|\xi_j|<2/m$, and let $f(x)=\int_{\bf R^n} \Phi(\xi)\delta(\xi_1+\cdots+\xi_n)e(-x\xi)d\xi$. Then for $h_j$ as in equation~\ref{eq2} we have,
\begin{align}
\sum h_1\left(\frac{\gamma_1}{T}\right) \ldots h_1\left(\frac{\gamma_n}{T}\right) f\left(\frac{L}{2\pi}\gamma_1,\ldots,\frac{L}{2\pi}\gamma_n\right) = \kappa({\bf h}) \frac{{\bf TL}}{2\pi} \int_{{\bf R^n}} C_{\underbar{o}}(v)\Phi(v)dv + O(T)
\end{align}
with
\begin{align*}
\int_{{\bf R^n}} \Phi(v)C_{\underbar{o}}(v)dv=\Phi(0)+\sum_{r=1}^{[n/2]} \sum \int\cdots\int|v_1|\ldots|v_r| \Phi(v_1e_{i(1)j(1)}+\cdots+v_re_{i(r)j(r)} dv_1\ldots dv_r
\end{align*}
where the sum over all choices of $r$ disjoint pairs of indices $i(t)<j(t)$ in \{1,\ldots,n\} and for $i<j$ we set
\begin{align*}
\begin{cases} e_{ij} = e_i-e_j \\
e_i=\{0,\ldots,1,0,\ldots\} {\text{the $i$-th standard basis vector}}
\end{cases}
\end{align*}
\end{theorem}
This proof is one of the main parts of Rudnick-Sarnak's result. They have given the proof in few steps and in each step there is a lemma. Here we give a brief overview of the proof along with the proof idea of those lemmas.

\medskip

{\bf Proof Idea:} They have started the proof with the Fourier transformation of the sum $C_n(f,h,T)$ and after that it has become
\begin{align*}
C_n(f,h,T)= \int_{{\bf R^n}} \prod_{j=1}^{n} \bigg\{\sum_{\gamma_j} h_j (\frac{\gamma_j}{T})e^{iL_{\gamma_j}\xi_j} \bigg\} \Phi(\xi) \delta(\xi_1+\cdots\xi_n)d\xi
\end{align*}
In the next step, they have changed the above equation into the sums over primes by using the Explicit formula and with the help of the test functions. The Explicit formula gives the the desired sums over primes representation of $C_n(f,h,T)$ and by expanding the products one can get that $C_n(f,h,T)$ is an alternating sum of terms of the form
\begin{align*}
C_{r,s} = \sum\frac{c(n_1)\ldots c(n_r)\overline{c(n_{r+1})}\ldots\overline{c(n_{r+s})}}{\sqrt{n_1\ldots n_{r+s}}}A_{r,s}(n,T)
\end{align*}
where $c(n)=\Lambda(n)a(n)$ and
\begin{align}\label{eq3}
A_{r,s}(n,T)=T^n\int\prod_{j=1}^{r} (T(L\xi_j+\log n_j))\prod_{j=r+1}^{r+s}g_j(T(L\xi_j-\log n_j)) \prod g_{j,T}(TL\xi_j)\Phi(\xi)\delta(\xi_1+\cdots\xi_n)d\xi
\end{align}
In expressing $C_n(f,h,T)$ as a sum various $C_{r,s}(T)$ they have got terms from all possible choices of $r$ of the factors to be $S_j^{+}(\xi_j)$ and $S_j^{-}(\xi_j)$ and the remaining $k=n-(r+s)$ factors to be $g_{j,T}(TL\xi_j)$.
Now we split the proof into the contribution of some lemmas. \\

{\bf First Lemma:} This lemma proves that the integrals defining $A_{r,s}(n,T)$ are rapidly convergent. More mathematically it states that
\begin{align*}
1. g_T(x)=\begin{cases}
\log T, |x| \ll \log\log T \\
\frac{1}{|x|}, |x|\gg \log\log T
\end{cases}
2. \int |g_T(x)|dx \ll \log T
\end{align*}
Proof of the both part are similar. By assuming $x>0$ and taking the expansion of the term $g_T(x)$ they have shifted the contour to the right. The fisrt gamma factor of the expansion of $g_T(x)$ is holomorphic for $Re(s)>0$. One can join the both gamma factors and using Striling formula $h(r)$ is rapidly decreasing and bounded, which proves the first lemma.\\

{\bf Second Lemma:} It says that if $\Phi$ is supproted in $|\xi_1+\ldots+\xi_n|\leq\frac{2-\delta}{m}$ then $A_{r,s}(n,T)=0$ unless $|n_j|\ll T$ and $n_1n_2\ldots n_{r+s}\ll T^{2-\delta}$. The proof is very simple, consider the integrand from Equation~\ref{eq3}, there exist an $\eta\in Supp\Phi$ so that $\sum_{j}\eta_j=0$ such that
\begin{align}\label{eq4}
\begin{cases}
|T(\eta_jL+\log n_j)| \ll 1, j=1,\ldots,r \\
|T(\eta_jL-\log n_j)|\ll 1, j+=r+1,\ldots,r+s
\end{cases}
\end{align}
%Hence, $\frac{n_j}{T^{m|\eta_j|}=1+O(\frac{1}{T})$ so that $n_j\ll T^{m|\eta_j|\ll T^{1-\delta/2}$ and $n_1n_2\cdotsn_{r+s}\ll T^{m\sum|\eta_j|}\ll T^{2-\delta}$. \\
Hence, $\frac{n_j}{T^{m|\eta_j|}}=1+O(\frac{1}{T})$ so that $n_j\ll T^{m|\eta_j|}\ll T^{1-\delta/2}$ and $n_1n_2\cdots n_{r+s}\ll T^{m\sum|\eta_j|}\ll T^{2-\delta}$. \\

{\bf Third Lemma:} This lemma tells about the Equation~\ref{eq3} for $j>r+s$ with restriction to $|TL\xi_j|\ll T^{\delta/3}$. Then for sufficiently large $T$ $\tilde{A}_{r,s}(n,T)$ is either 0 or it satisfies the second lemma. Also, $n_1n_2\ldots n_r=n_{r+1}\ldots n_{r+s}$. The proof is tricky using contradiction. Since $g$ is compactly supported, in order that the integrand not vanished consider some $\eta\in Supp \Phi$ satisifies equation~\ref{eq4}, with $|TL\eta_j|\ll T^{\delta/3}$ for $j>r+s$. Then the previous asumptions and some simple calculations proves the lemma. Actually in this lemma the target is to replace $A_{r,s}(n,T)$ by $\tilde{A}_{r,s}(n,T)$ to get the corresponding sum $\tilde{C}_{r,s}(n,T)$ in place of $C_{r,s}(n,T)$.\\

{\bf Fourth Lemma:} This lemma is very crucial. It calculates the region of integration for the difference $A_{r,s}(n,T)$ and $\tilde{A}_{r,s}(n,T)$ for $NM\ll T^{2-\delta}$. The region of the difference is a union of the sets $\mathcal{F}_{||}=\{\xi: |\mathcal{TL}|\xi\gg\mathcal{T}^{\delta/\ni}\}, k>r+s$. Observe that one can estimate the integral over the region $\mathcal{F}$. For this purpose, set
\begin{align*}
x_j=\begin{cases}
T(L\xi_j+\log n_j), 1\leq j\leq r \\
T(L\xi_j-\log n_j), r+1\leq j\leq r+s \\
TL\xi,
\end{cases}
\end{align*}
Using change of variable, previous lemmas and after some calculations by dividing the integral in $I_1$ and in $I_2$ one can conclude that

\begin{align}\label{eq5}
A_{r,s}(n,T)-\Tilde{A}_{r,s}(n,T)&=T^n\int_U\prod_{j=1}^r g_j(T(L\xi_j+\log n_j))\prod_{j=r+1}^{r+s}g_j(T(L\xi_j-\log n_j))\\ \nonumber
&\;.\prod_{j>r+s} g_{j,T}(TL\xi_j).\Phi(\xi)\delta(\xi_1+...+\xi_n)d\xi\\ \nonumber
&\ll\frac{T}{L^{n-1}}\int_{\{|x_j|\ll TL,\; |x_n|\gg T^{\delta/3}\}}\prod_{j\leq r+s}|g_j(x_j)|\prod_{r+s<j\leq n-1}|g_{j,T}(x_j)|\\ \nonumber
&\;.\left|g_{n,T}\left( T\log M-T\log N-\sum_{j=1}^{n-1}x_j\right)\right|dx_1...dx_{n-1}
\end{align}
After lots of computation by splitting the integral into $I_1$ and $I_2$ over region $|\sum x_j|\gg T|\log M/N|$. Computing and combining the integral $I_1$ and $I_2$ we can conclude this lemma. \\

Now we move our focus to compute the difference $C_{r,s}(T)-\tilde{C}_{r,s}(T)$, by dividing into diagonal and off-diagonal sum over $n$, which also implies that $M=N$, where $M=n_1\ldots n_r$ and $N=n_{r+1}\ldots n_{r+s}$. The diagonal sum over $n$ says that such that $|\log M/N|\ll T^{\delta/3}$ and the off-diagonal sum over $n$ says that $T|\log M/N|\gg T^{\delta/3}$. Computations says that $\sum_{diag}\ll T^{1-\delta/3}$ and $\sum_{off}\ll \frac{1}{L^{r+s}}\sum_{MN\ll T^{2-\delta},M\ne N} \frac{a_r(M)a_s(N)}{\sqrt{MN}|\log M/N|}$. \\
Our next focus is on the $a_k(m)$ sum, next lemma tells about it.\\

{\bf Fifth Lemma:} For $k\geq1$ fiixed, and any $\epsilon>0$
\begin{align*}
\sum_{m\leq X} a_k(m)^2\ll_{\epsilon} X^{1+\epsilon}.
\end{align*}
To prove the
above equation Rudnick-Sarnak has used two facts Chauchy-Schwartz and the number of ways of writing $m=m_1\cdots m_k$ is $O(m^{\epsilon})$ for any $\epsilon>0$. These facts help us to conclude that
\begin{align*}
\sum_{m\leq X} a_k(m)^2 \ll X^{\epsilon} \sum_{m_1\cdots m_k\leq X} |c(m_1)|^2\ldots|c(m_k)|^2.
\end{align*}
with $\sum_{n\leq X}|c(n)|^2\ll X^{1+\epsilon}$, for all $\epsilon>0$ follows from the absolute convergence of the differentiation of the partial L-function in $Re(s)>1$. By dyadic decomposition of the above equation they have concluded the result. From the diagonal term the main term contribution will come and the off-diangonal calculate the error term. The error term is not so sharp but it can give a good approximation of the error. Likewise, Rudnick-Sarnak have calculated the sum $\mathcal{I}_{\infty}$ and $\mathcal{I}_{\epsilon}$ where the main term is contributed by the first sum $\mathcal{I}_{\infty}$ and the second term is the error term. \\

{\bf Sixth Lemma:} This lemma approximate the function $g_i(y_j)$. Set, $k=n-(r+s)$, we have
\begin{align*}
\int_{V} \prod_{j=1}^{r+s} g_j(y_j) \prod_{j>r+s} g_{j,T}(y_j)dy=\frac{1}{2\pi} \kappa(h)L^k+O(L^{k-1})
\end{align*}
One can prove it by spliting the integral into two terms, first case is over $V$ and the second one is integrating over $\sum y_j=0$. Using Parseval inequality and Striling formula one can get the desired result. \\

{\bf Seventh and eighth Lemma:} This lemma calculates the sum $C_{r,s}(T)$ for $r+s>0$. For $r+s\geq3$, it proves that the sum over all indices of the function $c(p^r)$ is finite. It takes the Rankin-Selberg integration technique to prove this result. This lemma deals with sums over several prime factors. \\

{\bf Ninth Lemma:} It states that, if $1\leq r\leq s$ then
\begin{align*}
\sum_{\substack{p_1^{k_1}\ldots p_r^{k_r}=q_1^{l_1}\ldots q_s^{l_s} \\ p_i^{k_i},q_j^{l_j}\leq x}} \frac{c(p_1^{k_1})\cdots c(p_r^{k_r})\overline{c(q_1^{l_1})}\cdots\overline{c(q_s^{l_s})}}{p_1^{k_1}\cdots p_r^{k_r}}=\begin{cases}
O((\log x)^{2r}), r=s\\
O((\log x)^{2r-2}),
\end{cases}
\end{align*}
The proof idea is to divide the sum first into the sub-sums according to the number of distinct prime factors appearing and to collect the factors together which are corresponding to the same prime. Observe that, by doing this step the sum actually has become a sum of products which run over the primes with $p_i\ne p_j$. Clearly there are at most $r$ factors, and each factors contribute a bounded quantity unless $a=b=1$, which is the first case of this lemma. For the second case the product is $O(\log^{2r-2}x)$ unless $r=s$ and $a=b=1$, hence prove the second part. For the case $r=s$ each product contribute an error term and after multiplying them it comes a new error term which is the second case. Finally, the sum $C_{r,s}(T)=O(T)$ unless $r=s>0$. For another case $r=s$ next lemma state and prove the result. \\

{\bf Tenth Lemma:} It states that
\begin{align*}
&1. C_{r,s}(T)=O(T) \text{ unless } r=s>0. \\
&2. \text{ If $r=s$ then } C_{r,r}=\frac{TL}{2\pi} \kappa(h)\cdot\sum_{\sigma\in S_r} \int_{0}^{1/m} \cdots \int_{0}^{1/m} v_1\cdots v_r\cdot \\
&\Phi(-v_1,\dots, -v_r,v_{\sigma(1)},\ldots,v_{\sigma(r)},0\ldots,0)dv_1\ldots dv_r+O(T)
\end{align*}
where $S_r$ is the permutation group on $r$ letters.
The proof of the second equation can be done by summation by two parts. The distinct prime factors are going to contribute the error term, not the main term. The main term will be contributed by the permutation of distinct prime factors denoted by $q_j$. Therefore the condition of summing over distinct prime factors can be omitted.

\medskip

For $m=1$ or for the case of $\zeta(s)$ a polar term at $s=1$ will come and its corresponding coefficients $\lambda(n)$ in the sum are non-negative and the polar terms can not be dominated. Also, the central diagonal terms will come. In that case, Dirichlet L-functions $L(s,\chi)$, the same is true but no longer because it has no polar term. As we have said before the diagonal term will dominate the off-diagonal terms, so the main term will be contributed by the diagonal term. For the Riemann zeta function $\zeta(s)$ it is an exceptional case.

Now it's time to state another important theorem reagrding the asymptotic behaviour of $C_n(f,T)$.
\begin{theorem}
Let $\Phi\in C^2({\bf R^n})$ be supported in $\sum|\xi_j|<2/m$, and $f$ be given by
\begin{align*}
f(x)=\int_{{\bf R^n}} \Phi(\xi)\delta(\xi_1+\cdots+\xi_n)e(-x\cdot \xi)d\xi.
\end{align*}
Assume Riemann Hypothesis for $L(s,\pi)$ then
\begin{align*}
C_n(f,T)\tilde N(T)\int_{{\bf R^n}} \Phi(u) C_{\underline{o}}(u) du+O(T)
\end{align*}
\end{theorem}

%{\bf Proof Idea:} In theorem~\ref{thm2} take $h_j=\chi_{[-1,1]}$ be the characteristic function in the inerval $[-1,1]$. Assuming Riemann Hypothesis for $L(s,\pi)$ one can show that for the functions $h$ and $f$ as defined $C_n(f,T)$ is going to $\kappa(h)\mu(f)$, where $N(T)$ is the set of all $\gamma_j$'s less than equal to $T$, which of size $TL/2\pi$.

After obtatining the asymptotic behaviour of the $C_n(f,T)$ term Rudnick-Sarnak used the combinatorial sieve argumets to get the asymptotic behaviour of the $n$correlation sum $R_n(B_N,f)$ and they have concluded that
\[
R_n(B_N,f) \rightarrow \int_{{\bf R^n}} f(x) W_n(x)\delta\left(\frac{x_1+\cdots x_n}{n}\right)dx_1\ldots dx_n.
\]
as $N\to\infty$.

\section*{Acknowledgement(s)}

This work was supported by our guide and mentor Prof. Ritabrata Munshi. We thank our guide from the Tata Institute of Fundamental Research, Bombay (presently at Indian Statistical Institute, Kolkata) who provided insight and expertise that greatly assisted the research, although he may not agree with all of the interpretations of this article.
We would also like to show our gratitude to the Prof. Ritabrata Munshi, (TIFR, Bombay) for sharing his pearls of wisdom with us during the course of this work.

\end{document}